\documentstyle[12pt]{article}

  \newcommand{\const}{\rm const}
  \newcommand{\Var}{\rm Var}

  \newcommand{\card}{\rm card}

  \newcommand{\Sub}{\rm Sub}
  \newcommand{\StSub}{\rm StSub}

\begin{document}

   \begin{center}

 {\bf Exponential exact estimation for maximum  and minimum }\\

\vspace{4mm}

{\bf tail of distribution for non - Gaussian random vector.} \par

\vspace{6mm}

{\bf Formica M.R., Ostrovsky E., and Sirota L. } \par

\vspace{4mm}

\end{center}

 Universit\`{a} degli Studi di Napoli Parthenope, via Generale Parisi 13, Palazzo Pacanowsky, 80132,
Napoli, Italy. \\

e-mail: mara.formica@uniparthenope.it \\

 \vspace{3mm}

Department of Mathematics and Statistics, Bar-Ilan University,
59200, Ramat Gan, Israel. \\

\vspace{3mm}

e-mail:eugostrovsky@list.ru\\
Department of Mathematics and Statistics, Bar-Ilan University,\\
59200, Ramat Gan, Israel.

\vspace{3mm}

e-mail:sirota3@bezeqint.net \\

\vspace{4mm}

\begin{center}

  {\bf Abstract} \par

\vspace{3mm}

 \end{center}

 \hspace{3mm} We find the exponential exact two -terms non - asymptotic expression for the maximum  and minimum distribution
 of a non - Gaussian, in general case, random vector.

\vspace{3mm}

\begin{center}

 \ {\it Key words and phrases:}

 \vspace{3mm}

\end{center}

  \hspace{3mm} Random vector, distribution of maximum, Lebesgue - Riesz and Grand Lebesgue Spaces and norms,
  exact exponential tail estimate, Bonferroni's inequality, H\"older's inequality, Young - Fenchel transform,
  moment generation function,  subgaussian variables, vectors and correlation;  triangle inequality. \\

 \vspace{5mm}

 \section{Statement of problem. Notations. Foreword. }

  \vspace{5mm}

\hspace{3mm} Let $ \ (\Omega = \{\omega\}, \cal{B}, {\bf P})  \ $ be certain (sufficiently rich) probability space and let
$ \  D = \{1,2,\ldots,d\} \ $ be a positive integer numerical set,
$ \ \xi = \vec{\xi} = \{ \xi_1, \xi_2, \ldots, \xi_d \}, \ d = 2,3,\ldots \ $ be a random vector. Define the following important notion of the
tail of maximum distribution

\begin{equation} \label{particular max}
Q[S](u) \stackrel{def}{=} {\bf P} \left[ \max_{i \in S} \xi_i > u \ \right],
\end{equation}
so that

\vspace{3mm}

\begin{equation} \label{Important notion}
Q(u) = Q[D](u) \stackrel{def}{=} {\bf P} \left[ \max_{i = 1,2,\ldots,d}\xi_i > u \right]
\end{equation}
for the sufficiently great values $ \ u, \ $ say for $ \ u \ge 1. \ $ \par
 \ In particular,

$$
Q_i(u) := {\bf P}(\xi_i > u), \hspace{3mm} Q_{i,j}(u) := {\bf P}(\xi_i > u, \ \xi_j > u), \ i,j \in D, \ i \ne j.
$$

\vspace{3mm}

 \hspace{3mm} The case of the minimum distribution

 \begin{equation} \label{def min  tail}
 S[\xi](u) \stackrel{def}{=} {\bf P} ( \ \min_{i \in D} \xi_i > u  \ ), \ u \ge 1.
 \end{equation}

 will be considered further.

 \vspace{3mm}

 \hspace{3mm}  {\bf  Our purpose in this report is to deduce the  exact exponential tail estimations for these tail probabilities. }\par

 \vspace{3mm}

  \ {\it We do not suppose the Gaussian distribution of the source vector} $ \ \xi. \ $ \\

\vspace{4mm}

 \ There are a huge numbers of works devoted to this problem, see e.g.
\cite{anatrielloformicaricmat2016}, \cite{Belyaev-Piterbarg1978},  \cite{Buld Koz AMS},  \cite{Ermakov etc. 1986},
\cite{fioguptajainstudiamath2008},  \cite{Ng-Tang-Yang}, \cite{Ostrovsky1999},
\cite{Ostrovsky 1994},  \cite{Pisier 1980},  \cite{Sgibnev 1996} etc. The applications of these estimates in the Law of Iterated
Logarithm (LIL) and following in statistics are investigated in particular  in \cite{Hu 1991}, \cite{Sgibnev 1996}, \cite{Stout 1974},
\cite{Sung 1996}.\par
 \ Note that as a rule in the mentioned works is considered the case when $ \ d \to \infty. \ $ \par

\begin{center}

 \vspace{4mm}

 \hspace{3mm} {\sc Preliminary estimations.} \par

 \vspace{4mm}

 \end{center}

 \hspace{3mm}  We will use the classical Bonferroni's inequality

\begin{equation} \label{key ineqalities}
\sum Q_i(u) - \sum \sum Q_{i,j}(u) \le Q[D](u) \le \sum Q_i(u), \ u \ge 1.
\end{equation}

\vspace{4mm}

  where by definition

 $$
   \sum \stackrel{def}{=} \sum_{i \in D}, \hspace{3mm}  \sum \sum \stackrel{def}{=} \sum \sum_{i,j \in D; \ i \ne j}.
 $$

 \vspace{5mm}

 \section{Grand Lebesgue Spaces approach.}

 \vspace{5mm}

 \hspace{3mm}  We present here for beginning some known facts from the theory of one - dimensional random variables with exponential decreasing tails of distributions, see
\cite{Buld Koz AMS}, \cite{Buldygin-Mushtary-Ostrovsky-Pushalsky}, \cite{caponeformicagiovanonlanal2013}, \cite{Ermakov etc. 1986},
\cite{Fiorenza-Formica-Gogatishvili-DEA2018}, \cite{fioforgogakoparakoNAtoappear}, \cite{Iwaniec-Sbordone 1992}, \cite{Kozachenko-Ostrovsky 1985},
\cite{Kozachenko at all 2018}, \cite{liflyandostrovskysirotaturkish2010}, \cite{Ostrovsky1999}, \cite{Ostrovsky 1994}, \cite{Ostrovsky HIAT}. \par

 \ Let $ \ \phi = \phi(\lambda), \ \lambda \in (-\lambda_0, \lambda_0), \ \exists \lambda_0 = \const \in (0,\infty] \ $
be certain even strong convex which takes positive values for positive arguments twice continuous differentiable
function, briefly: Young-Orlicz function, such that

$$
\phi(0) = 0, \ \phi'(0) = 0, \ \phi^{''}(0) \in (0,\infty).
$$

 \ For instance: $ \ \phi(\lambda) = 0.5 \ \lambda^2, \ \lambda \in R; \ - \ $ the so-called subgaussian case. \par

\ We denote the set of all these Young-Orlicz function as  $ \ \Phi: \ \Phi = \{ \phi(\cdot)  \}. $
We say by definition that the centered random variable (r.v) $ \ \xi \ $
belongs to the space  $ B(\phi), \ $  if there exists certain non-negative constant \ $ \tau \in [0,\ \infty), \ $ such that

\begin{equation} \label{Bphi set}
\forall \lambda \in (-\lambda_0, \lambda_0) \ \Rightarrow {\bf E} \exp(\pm \lambda \ \xi) \le \exp(\phi(\lambda \ \tau)).
\end{equation}

\vspace{5mm}

\ {\bf Definition 2.1.} \ The minimal non-negative value $ \ \tau \ $ satisfying the last relation (\ref{Bphi set})
 for all the values $ \ \lambda \in (-\lambda_0, \ \lambda_0), \ $ is named as a $ \ B(\phi) \  $  norm of the variable $ \ \xi, \ $  write
$  \ ||\xi||B(\phi) \stackrel{def}{=} \ $

\begin{equation} \label{important definition}
\inf \{\tau, \ \tau \ge 0, \  \forall \lambda \in (-\lambda_0, \lambda_0) \ \Rightarrow {\bf E} \exp(\pm \lambda \ \xi) \le \exp(\phi(\lambda \ \tau)) \ \} .
\end{equation}

\vspace{4mm}

\ For instance: $ \ \phi(\lambda) = \phi_2(\lambda) : = 0.5 \ \lambda^2, \ \lambda \in R; \ - \ $ the so-called subgaussian case, \ write $ \ \xi \in \Sub. \ $ \par
 \ We will denote as ordinary the norm in this very popular space by $ \  ||\cdot||\Sub: \ $

$$
||\xi||\Sub \stackrel{def}{=} ||\xi||B(\phi_2).
$$
 \ It is known  that if the r.v. $ \ \xi_i \ $ are independent and subgaussian, then

$$
||\sum_{i=1}^n \xi_i||\Sub \le \sqrt{ \sum_{i=1}^n ||\xi_i||^2 \Sub}.
$$

\vspace{4mm}

 \ {\bf Definition 2.2.} The centered r.v. $ \ \eta \ $ with finite non - zero variance $ \  \sigma^2 := \Var (\eta) \in (0,\infty)  \  $ is said to be strictly subgaussian,
write: $ \ \eta \in \StSub, \ $ iff

\begin{equation} \label{StSub}
{\bf E} \exp(\pm \ \lambda \ \eta) \le \exp(0.5 \ \sigma^2 \ \lambda^2), \ \lambda \in R.
\end{equation}

 \ For instance, every centered non - zero Gaussian r.v. belongs to the space $ \ \StSub. \ $ The Rademacher's r.v. $ \ \eta: \ {\bf P}(\eta = 1) = {\bf P}(\eta = - 1) = 1/2  \ $
is also strictly subgaussian. Many other strictly subgaussian r.v. are represented in \cite{Buld Koz AMS},  \cite{Ostrovsky1999}. \par

\vspace{4mm}

 \ It is known that the set $ \   B(\phi), \ \phi  \in \Phi \ $ relative the norm (2.3) and ordinary algebraic operations forms a Banach functional
rearrangement invariant space, which are equivalent the so - called Grand Lebesgue ones as well as Orlicz exponential spaces.
These spaces are very convenient for the investigation of the r.v. having an exponential decreasing tail of distribution, for instance, for investigation of the limit
theorem, the exponential bounds of distribution for sums of random variables, non-asymptotical properties, problem of continuous and weak compactness of random
fields, study of Central Limit Theorem in the Banach space etc. \par

 \ Denote by $ \ \nu(\cdot)\ $ the Young - Fenchel, or Legendre transform for the function $ \ \phi: \ $

\begin{equation} \label{Young Fenchel}
\nu(x) = \nu[\phi](x)  \stackrel{def}{=} \sup_{\lambda|: |\lambda\le \lambda_0} (\lambda x - \phi(\lambda)) = \phi^*(x).
\end{equation}

\vspace{5mm}

 \ It is important for us in this  preprint to note that if the non - zero r.v. $ \ \xi \ $ belongs  to the space $ \ B(\phi) \ $ then

\begin{equation} \label{tail Bphi}
{\bf P}(\xi > x) \le \exp \left(  - \nu(x/||\xi||B(\phi)) \  \right).
\end{equation}
 \ The inverse conclusion is also true up to multiplicative constant under  suitable conditions.\par

\vspace{5mm}

 \ Further, assume that the {\it centered} r.v. $ \ \xi \ $  has a finite in some non - trivial neighborhood of origin finite moment generation function

$$
\phi[\xi](\lambda) \stackrel{def}{=} \max_{\pm} \ \ln {\bf E} \exp ( \ \pm \lambda \ \xi \ ) < \infty, \ \lambda \in ( - \ \lambda_0, \lambda_0)
$$
for some $ \ \lambda_0 = \const \in (0, \ \infty]. \ $  Obviously, the last condition is quite equivalent to the well - known Cramer's one. \par

 \ We can and will agree $ \  \phi[\xi](\lambda) := \infty  \  $ for all the values $ \ \lambda \ $ for which

$$
 {\bf E} \exp ( \ |\lambda| \ \xi) = \infty.
$$

 \ The introduced in (2.6) function $ \ \phi_{\xi}(\lambda) \ $ is named as {\it  natural } function for the r.v. $ \ \xi; \ $ herewith  $ \ \xi \in B(\phi[\xi]) \ $ and moreover

$$
||\xi||B(\phi[\xi]) = 1.
$$

\vspace{5mm}

\section{Main result. Upper maximum tail estimate.}

\vspace{5mm}

 \hspace{3mm} Let us return to the formulated  above problem, $ \ Q[\xi](u) \ $ estimation. We suppose that there exists certain function
 $ \  \phi \in \Phi  \ $ such that

 \begin{equation} \label{betai}
 \forall i \in D \ \Rightarrow \ \beta_i := \ || \xi_i||B\phi \in (0,\infty).
 \end{equation}

 \ Set also

\begin{equation} \label{beta r}
\beta:= \max_{i \in D}  \beta_i, \ r := \card \{ \ \beta_i: \ \beta_i = \beta \ \},
\end{equation}
where as ordinary $ \ \card(M) \ $ denotes the amount  of elements of the set $ \ M; \ $ \par

$$
R = R[\beta] \stackrel{def}{=} \{i; \ i \in D, \ \beta_i < \beta \},
$$

$$
\underline{\beta} \stackrel{def}{=} \max \{\beta_i, \beta_i < \beta\} = \max \{\beta_i, \ i \in R[\beta] \}.
$$

\vspace{4mm}

 \ {\bf Proposition  3.1. } We conclude by means of the Bonferroni's inequality  (\ref{key ineqalities})

 \begin{equation} \label{upper est 1}
 Q[\xi](u) \le r \ \exp( \ - \nu[\phi](u/\beta) \ ) + \sum_{i \in R} \exp(-\nu[\phi](u/\beta_i)),
 \end{equation}
and following

 \begin{equation} \label{upper est 2}
 Q[\xi](u) \le r \ \exp( \ - \nu[\phi](u/\beta) \ ) + (d - r) \exp(-\nu[\phi](u/\underline{\beta})).
 \end{equation}

 \vspace{4mm}

 \section{Main result: lower maximum tail estimates.}

\vspace{4mm}

 \hspace{3mm} We intent  to apply  here  the  left - hand side of  the
 relation of Bonferroni (\ref{key ineqalities}) and we retain the  restrictions
(\ref{betai}). Moreover, let us suppose

\begin{equation} \label{delta cond}
\exists \delta_i \in (0,\beta_i] \ \Rightarrow Q_i(u) \ge \exp \left\{ \ - \nu[\phi](u/\delta_i)  \ \right \}, \ u \ge 1.
\end{equation}

 \vspace{3mm}

 \ Further, we deduce applying the triangle inequality for the $ \ B(\phi) \ $ norm

 $$
  ||\xi_i + \xi_j||B\phi \le \beta_i +\beta_j, \ i \ne j,
 $$
that

\begin{equation} \label{double ij}
Q_{i,j}(u) \le {\bf P}(\xi_i + \xi_j \ge 2 u) \le \exp \left\{ \ - \nu[\phi](2u/(\beta_i + \beta_j)) \ \right\}.
\end{equation}

 \vspace{4mm}

 \hspace{3mm} To summarize. \par

 \vspace{3mm}

 \ {\bf Proposition 4.1.}

 \vspace{3mm}

 \begin{equation} \label{beginning lower}
  \ Q[\xi](u) \ge \  \sum \ \exp \left\{ \ - \nu[\phi](u/\delta_i)  \ \right \} -
 \end{equation}

\begin{equation} \label{lower main}
  \sum \ \sum \ \exp \left\{ \ - \nu[\phi](2u/(\beta_i + \beta_j)) \ \right\}, u \ge 1.
\end{equation}

 \vspace{5mm}

 \section{ A note about the minimum distribution.}

\vspace{5mm}

 \hspace{3mm} Let us take into account the tail of the {\it minimum  distribution} for the source  random vector (we recall)

 \begin{equation} \label{def min distr}
 S[\xi](u) \stackrel{def}{=} {\bf P} ( \ \min_{i \in D} \xi_i > u  \ ), \ u \ge 1.
 \end{equation}

\vspace{3mm}

 \ We assume that the r.v. $ \  \xi = \vec{\xi} \ $ is non - negative: $ \ \xi_i \ge 0 \ $ and  define the so - called
 multivariate moment generation function (cf. with a characteristical function)

\begin{equation} \label{MGF}
g[\xi](\lambda) = g(\lambda) \stackrel{def}{=} {\bf E} \exp(\lambda,\xi),
\end{equation}
where $ \ \lambda = \vec{\lambda} \in R^d, \ (\lambda,\xi) = \sum \lambda_i \xi_i; \ $
and as ordinary $ \ (\lambda,\xi) = \sum \lambda_i \xi_i, \ |\lambda| := \sqrt{(\lambda,\lambda)}.  \ $ \par
 \ It will be presumed that this function is finite at last in some neighborhood of origin:

\begin{equation} \label{finiteness}
\exists \epsilon \in (0, \infty]  \  \forall \lambda: \ |\lambda| < \epsilon \ \Rightarrow g(\lambda) < \infty.
\end{equation}

 \ Of course, if the random vectors $ \ \xi,\eta  \ $ are independent, then $ \  g[\xi + \eta](\lambda) = g[\xi](\lambda) \cdot g[\eta](\lambda).\ $\par

\vspace{3mm}

 \ Let now $ \ q_i, \ i \in D \ $  be arbitrary numerical vector such that

$$
q_i \ge 1, \ \sum \frac{1}{q_i} = 1.
$$
 \ The H\"older's inequality give us

$$
g[\xi](\lambda) \le \prod_{i \in D} \left[ \ {\bf E} \exp(  \ q_i \ \lambda_i \ \xi_i) \ \right]^{1/q_i} = \prod_{i \in D}  \left\{ \ g[\xi_i](\lambda_i \ q_i) \ \right\}^{1/q_i}.
$$

 \ In particular,

$$
g[\xi](\lambda) \le \left\{ \ \prod_{i \in D} g[\xi_i](\lambda_i \ d) \ \right\}^{1/d}.
$$

 \ If in addition the r.v. $ \ \{\xi_i\} \ $ are identical distributed  and $ \ \lambda_i = \mu = \const > 0, \ $ then

$$
g[\xi](\lambda) \le g[\xi_1](\mu \ d).
$$

\vspace{4mm}

 \ Further,  define as usually the so - called Young - Fenchel, or Legendre transform in the multivariate case
 $ \  x = \vec{x} = \{x_i\} \in R^d \ $

 \vspace{3mm}

 \begin{equation} \label{Young Fen}
 g^*(\vec{x}) = g^*(x) \stackrel{def}{=} \sup_{\lambda \in R^d} [( \lambda,x) - g(\lambda)].
 \end{equation}

\vspace{3mm}
 It is known, see \cite{Buld Koz AMS}, \cite{Ostrov Prokhorov} that under formulated restrictions

\begin{equation} \label{tail vec estim  beg}
{\bf P} (\forall i \in D \ \Rightarrow  \ \xi_i \ge x_i)   = {\bf P} \left[ \ \cap_{i \in D} \{ \xi_i \ge x_i  \}  \ \right] \le
\end{equation}

\begin{equation} \label{tail vec estim  en}
\exp(-g[\xi]^*(\vec{x})), \ x_i \ge 0.
\end{equation}

\vspace{4mm}

 As a consequence: \par

 \ {\bf Proposition 5.1.}

\begin{equation} \label{min estim}
S[\xi](u) \le \exp \left\{ \ - g^*[\xi](u,u, \ldots,u)  \ \right\}, \ u \ge 1.
\end{equation}

\vspace{5mm}

\section{Examples.}

\vspace{5mm}

\ {\bf Definition 6.1.} The centered two - dimensional random vector $ \  \xi = \vec{\xi} = (\xi_1, \xi_2) \  $ is said to be
{\it subgaussian}, if

\begin{equation} \label{two dim sub}
{\bf E} \exp(\lambda_1 \ \xi_1 + \lambda_2 \ \xi_2) \le\exp \left\{ \ 0.5(\sigma_1^2 \lambda_1^2 + 2 \rho \ \sigma_1 \ \sigma_2 \ \lambda_1  \ \lambda_2 + \sigma_2^2 \lambda_2^2) \ \right\}
\end{equation}
for arbitrary real numbers $ \ \lambda_1,\lambda_2 \in R. \ $  Here $ \ \sigma_1, \sigma_2 = \const > 0, \  \rho = \const, \ |\rho| < 1.  \ $ \par
 \ For instance, the relation (\ref{two dim sub}) is satisfied for arbitrary non - degenerate, i.e. having strictly positive  definite (symmetrical) covariation
 matrix,  centered two - dimensional Gaussian distributed
 random vector, as well as for the independent subgaussian centered random variables $ \  (\xi_1, \xi_2).  \  $ The value $ \ \rho \ $ may be named as a
 {\it subgaussian correlation}  between r.v. $ \ \xi_1, \ \xi_2. \ $ \par
 \ We conclude by virtue of Proposition 5.1 after simple calculations for the values in particular  $ \ u \ge 1 \ $

\begin{equation} \label{min two dim}
{\bf P}(\min(\xi_1, \xi_2) > u) \le \zeta[\sigma_1,\sigma_2,\rho](u),
\end{equation}
where

\begin{equation} \label{zeta}
\zeta[\sigma_1,\sigma_2,\rho](u) \stackrel{def}{=}
\exp \left\{ \ - \frac{u^2}{2 (1 - \rho^2)} \ \frac{\sigma_1^2 + \sigma_2^2 - 2 \rho \ \sigma_1 \ \sigma_2}{\sigma_1^2 \sigma_2^2} \ \right\}.
\end{equation}

\vspace{3mm}

 \ Of course,

\begin{equation} \label{ultimately est}
{\bf P}(\min(\xi_1, \xi_2) > u) \le \min \{ \ \zeta[\sigma_1,\sigma_2,\rho](u), \ \exp \left\{ \ - g^*[\xi_1,\xi_2](u,u)  \ \right\} \ \}.
\end{equation}

\vspace{3mm}

 \ {\bf Remark 5.1.} The last estimate  (\ref{ultimately est}) is  exponential asymptotically as $ \ u \to \infty \ $ exact for instance for the
 centered gaussian distributed two - dimensional random vector. \par

\vspace{3mm}

 \ {\bf Remark 5.2.} This estimate is also asymptotically as $ \ u \to \infty \ $  exponential exact for the subgaussian r.v. in both the extremal
 cases $ \ \rho = 0 \ $ and $ \ \rho = 1 - 0. \ $ Indeed, let $ \ \xi_1, \xi_2 \ $ be  {\it independent} and standard subgaussian:

 \vspace{3mm}
$$
{\bf E} \xi_{1,2} = 0, \  ||\xi_1||\Sub = ||\xi_2||\Sub = 1.
$$
 \vspace{3mm}

 \ The relation (\ref{ultimately est} ) give us the following exponential exact estimate

 \vspace{3mm}

$$
 {\bf P}(\min(\xi_1, \xi_2) > u) \le \exp \left(   \  - u^2 \ \right), \ u \ge 1.
$$

 \vspace{4mm}

 \ Further, let now $ \  \xi_1 = \xi_2 = \xi \ $ be a standard subgaussian variable; then

$$
\min(\xi_1,\xi_2) = \xi, \ \rho = + 1.
$$

\ The relation (\ref{ultimately est} ) give us the following exponential exact estimate yet in this limiting case
$ \ \rho = 1 - 0 \ $

 \vspace{3mm}

$$
 {\bf P}(\min(\xi_1, \xi_2) > u) \le \exp \left(   \  - u^2/2 \ \right), \ u \ge 1.
$$

\vspace{3mm}

 \ {\bf Sub - remark 5.2.} Notice in addition to the foregoing remark 5.2 that in these conditions, aside from one
 that $ \ \rho = + 1 \ $ we assume instead $ \  \rho = -1; \ $ i.e. that $ \ \xi_1 = - \xi_2. \ $ This case is
 trivial for us, as long as under these restrictions

 $$
 {\bf P} ( \  \min(\xi_1,\xi_2) > u \ ) = 0, \ u > 0.
 $$

 \vspace{5mm}

\section{Concluding remarks.}

\vspace{5mm}

 \hspace{3mm} It is interest in our opinion to find the exact exponential {\it lower}  tail estimations for the minimum of
 random variables. \par

\vspace{6mm}

\vspace{0.5cm} \emph{Acknowledgement.} {\footnotesize The first
author has been partially supported by the Gruppo Nazionale per
l'Analisi Matematica, la Probabilit\`a e le loro Applicazioni
(GNAMPA) of the Istituto Nazionale di Alta Matematica (INdAM) and by
Universit\`a degli Studi di Napoli Parthenope through the project
\lq\lq sostegno alla Ricerca individuale\rq\rq .\par


\vspace{6mm}

\begin{thebibliography}{79}

\bibitem{anatriellofiojmaa2015}
{\bf G.~Anatriello} and {\bf A.~Fiorenza.} {\it Fully measurable
grand Lebesgue spaces}. J. Math. Anal. Appl. \textbf{422} (2015),
no.~2, 783--797.

\bibitem{anatrielloformicaricmat2016}
{\bf G.~Anatriello} and {\bf M.~R.~Formica.} {\it Weighted fully
measurable grand Lebesgue spaces and the maximal theorem}. Ric. Mat.
\textbf{65} (2016), no.~1, 221--233.

\bibitem{Belyaev-Piterbarg1978} {\bf Yu.~K.~Belyaev} and {\bf V.~I.~Piterbarg.} {\it Random processes. Sample paths and intersections}. Collection of articles, Publishing House "MIR",
Moscow (1978); 249--257, (in Russian).


\bibitem{Braverman1991} {\bf M.~Sh.~Braverman.} {\it  Bounds on the sums of independent random variables in symmetric spaces.} Ukrainian Mathematical Journal,
 \textbf{43} (1991), no.~2, 148--153.


\bibitem{Buld Koz AMS}
{\bf Buldygin V.V., Kozachenko Yu.V. }  {\it Metric Characterization of Random
 Variables and Random Processes.} 1998, \ Translations of Mathematics Monograph, AMS, v.188.



\bibitem{Buldygin-Mushtary-Ostrovsky-Pushalsky} {\bf V.~V.~Buldygin, D.~I.~Mushtary, E.~I.
~Ostrovsky} and {\bf M.~I.~Pushalsky.} {\it New Trends in
Probability Theory and Statistics.} Mokslas (1992), V.1, 78--92;
Amsterdam, Utrecht, New York, Tokyo.

\bibitem{caponeformicagiovanonlanal2013}
{\bf C.~Capone, M.~R.~Formica} and {\bf R.~Giova.} {\it Grand
{L}ebesgue spaces with respect to measurable functions}. Nonlinear
Anal. \textbf{85} (2013), 125--131.



\bibitem{Ermakov etc. 1986}
{\bf S. V. Ermakov, and E. I. Ostrovsky.} {\it Continuity Conditions, Exponential Estimates, and the Central Limit Theorem for Random Fields.}
 Moscow, VINITY,  1986. (in Russian).



\bibitem{Fiorenza2000} {\bf A.~Fiorenza.} {\it Duality and reflexivity in grand Lebesgue
spaces.} Collect. Math. \textbf{51} (2000), no. 2, 131--148.

\bibitem{fiokarazanalanwen2004}
{\bf A.~Fiorenza} and {\bf G.~E.~Karadzhov.} {\it Grand and small
Lebesgue spaces and their analogs}, Z. Anal. Anwendungen \textbf{23}
(2004), no.~4, 657--681.

\bibitem{fioguptajainstudiamath2008}
{\bf A.~Fiorenza, B.~Gupta} and {\bf P.~Jain.} {\it The maximal
theorem for weighted grand Lebesgue spaces}. Studia Math.
\textbf{188} (2008), no.~2, 123--133.


\bibitem{Fiorenza-Formica-Gogatishvili-DEA2018}
{\bf A.~Fiorenza, M.~R.~Formica} and {\bf A.~Gogatishvili.} {\it On
grand and small Lebesgue and Sobolev spaces and some applications to
PDE's}. \emph{Differ. Equ. Appl.} \textbf{10} (2018), no.~1, 21--46.

\bibitem{fioforgogakoparakoNAtoappear}
{\bf A.~Fiorenza, M. R.~Formica, A.~Gogatishvili, T.~Kopaliani} and
{\bf J.~M. Rakotoson.} {\it Characterization of interpolation
between grand, small or classical Lebesgue spaces}. Nonlinear Analysis, Vol.177, Part {\bf B,} Dezember 2018,
pages 422 \ - \ 453.


\bibitem{fioformicarakodie2017}
{\bf A.~Fiorenza, M.~R.~Formica} and {\bf J.~M. Rakotoson.} {\it
Pointwise estimates for {$ \ G\Gamma$c \ }-functions and applications}.
Differential Integral Equations, \textbf{30}, (2017), no. \ ~11 \ - \ 12,
809 \ - \ 824.

\bibitem{formicagiovamjom2015}
{\bf M.~R. Formica} and {\bf R.~Giova.} {\it Boyd indices in
generalized grand Lebesgue spaces and applications}. Mediterr. J.
Math., \textbf{12}, (2015), no.~3, 987 \ - \ 995.

\bibitem{Hu 1991}
{\bf T.~C.~Hu.} {\it  On the law of the iterated logarithm for
arrays of random variables}. Comm. Statist. Theory Methods
\textbf{20} (1991), no.~7, 1989--1994.

\bibitem{Hu-Weber 1992}
{\bf T.~C.~Hu} and {\bf N.~C.~Weber.} {\it On the rate of
convergence in the strong law of large numbers for arrays.} Bull.
Austral. Math. Soc., \textbf{45} (1992), no.~3, 479--482.

\bibitem{Iwaniec-Sbordone 1992}
{\bf T.~Iwaniec} and {\bf C.~Sbordone.} {\it On the integrability of
the Jacobian under minimal hypotheses.} Arch. Rational Mech. Anal.
\textbf{119} (1992), no.~2, 129--143.

\bibitem{Kolmogoroff 1929}
 {\bf A.~N.~Kolmogoroff.} {\it \"Uber das Gesetz des iterierten Logarithmus".} (German) Math. Ann., \textbf {101} (1929), no.~1, 126--135.

\bibitem{Kozachenko-Ostrovsky 1985}
{\bf Yu.~V.~Kozachenko} and {\bf E.~I.~Ostrovsky.} {\it The Banach
Spaces of random variables of sub-Gaussian type.} of Probab. and
Math. Stat., \textbf{32} (1985), (in Russian). Kiev, KSU, 43--57.


\bibitem{Kozachenko at all 2018}
{\bf Yu.V. Kozachenko, Yu.Yu. Mlavets, and N.V. Yurchenko.} {\it Weak convergence of stochastic processes from spaces} $F_\psi(\Omega).$
STATISTICS, OPTIMIZATION AND INFORMATION COMPUTING, Vol.6,  June 2018, pp. 266 - 277.


\bibitem{liflyandostrovskysirotaturkish2010}
{\bf E.~Liflyand, E.~Ostrovsky} and {\bf L.~Sirota.} {\it Structural
properties of bilateral grand {L}ebesque spaces}. Turkish J. Math.
\textbf{34} (2010), no.~2, 207--219.

\bibitem{Ng-Tang-Yang}
{\bf K.~W.~Ng, Q.~H.~Tang} and {\bf H.~Yang.} {\it Maxima of Sums of
Heavy-Tailed Random Variables.} ASTIN Bulletin: The Journal of the
IAA, {\textbf 32,} (2002), no.~1, pp. 43--55.


\bibitem{Ostrovsky1999}
{\bf E.Ostrovsky.} {\it Exponential estimates for random fields and
its applications.} 1999, OINPE, Moscow - Obninsk.

\bibitem{Ostrovsky 1994}
{\bf E.Ostrovsky.} {\it Exponential estimate in the Law of Iterated
Logarithm in Banach Space.} Math. Notes \textbf{56} (1994), no.~
5-6, 1165--1171.

\bibitem{Ostrovsky HIAT}
 {\bf Ostrovsky E. and Sirota L.} {\it Moment Banach spaces: theory and applications.}
HIAT Journal of Science and Engineering, C, Volume 4, Issues 1 - 2, pp. 233 \ - \ 262, (2007).


\bibitem{Ostrov Prokhorov}
{\bf  Ostrovsky E., Sirota L.} {\it Prokhorov-Skorokhod continuity of random fields.
A natural approach.} \\
arXiv:1710.05382v1 [math.PR] 15 Oct 2017













\bibitem{Pickands1967}
{\bf  J.~Pickands.} {\it  Maxima of stationary Gaussian processes.}
Z. Wahrscheinlichkeitstheorie und Verw. Gebiete \textbf{7} (1967),
190--223.

\bibitem{Pisier 1980}
{\bf G.~Pisier.} {\it Conditions d'entropie assurant la continuit\'e
de certains processus et applications \`{a} l'analyse harmonique.}
(French) Seminaire d'analyse fonctionnelle (1980), Exp. No. 13-14,
pp. 43 \ 46.

\bibitem{Qi 1994}
{\bf  Y.~C.~Qi.} {\it On strong convergence of arrays.} Bull.
Austral. Math. Soc. \textbf{50} (1994), no.~2, 219--223.


\bibitem{Samko-Umarkhadzhiev}
{\bf S.~G.~Samko} and {\bf S.~M.~Umarkhadzhiev.} {\it On
Iwaniec-Sbordone spaces on sets which may have infinite measure.}
Azerb. J. Math. \textbf{1} (1)\ (2011), 67--84.

\bibitem{Samko-Umarkhadzhiev-addendum}
{\bf S.~G.~Samko} and {\bf S.~M.~Umarkhadzhiev.} {\it On
Iwaniec-Sbordone spaces on sets which may have infinite measure:
addendum.} Azerb. J. Math \textbf{1} (2) \ (2011), 143--144.

\bibitem{Sgibnev 1996}
{\bf M.~S.~Sgibnev.} {\it On the distribution of the maxima of
partial sums.} Statist. Probab. Lett. \textbf{28} (1996), no.~3,
235--238.

\bibitem{Stout 1974}
{\bf W.~F.~Stout.} {\it Almost sure convergence.} Probability and
Mathematical Statistics, Vol. \textbf{24}. Academic Press. New
York-London, 1974.

\bibitem{Sung 1996}
{\bf  S.~H.~Sung.} {\it  An analogue of Kolmogorov's Law of the
Iterated  Logarithm for arrays.} Bull. Austral. Math. Soc.
\textbf{54} (1996),  no. 2, 177--182.




\bibitem{Teicher 1981}
{\bf  H.~Teicher.} {\it Almost certain behavior of row sums of
double arrays.} Analytical methods in probability theory
(Oberwolfach, 1980),  pp. 155--165, Lecture Notes in Math.,
\textbf{861}, Springer, Berlin-New York, 1981.

\end{thebibliography}
\end{document}